\magnification = 1150
\showboxdepth=0 \showboxbreadth=0

\baselineskip 14pt
\parskip3pt
\def\qed{\hfill\vrule height6pt width6pt depth0pt}

\def\ss{\smallskip}
\def\ms{\medskip}
\def\bs{\bigskip}

\def\cl{\centerline}

\def\nind{\noindent}

\def\ref#1#2{\nind\hangindent.5in\hbox to .5in{#1\hfill}#2}
\def\reff#1#2{\nind\hangindent.8in\hbox to .8in{\bf #1\hfill}#2\par}
\def\refd#1#2{\nind\hangindent.8in\hbox to .8in{\bf #1\hfill {\rm
--}}#2\par}
\def\pmb#1{\setbox0=\hbox{#1}
\kern-0.025em\copy0\kern-\wd0
\kern.05em\copy0\kern-\wd0
\kern-.025em\raise.0433em\box0}
\def\ca#1{{\cal #1}}
\def\bo#1{{\bf #1}}

\def\Cal{\ca}     

\def\R{\bo R}

\def\frac#1#2{{#1\over#2}}

\def\text#1{\rm#1}

\font\au=cmcsc10

\outer\def\stmnt  #1. #2\par{\medbreak
\noindent{\bf#1.\enspace}{\sl#2}\par
\ifdim\lastskip<\medskipamount \removelastskip\penalty55\medskip\fi}

\def\newpage{\vfill\eject}
\def\newline{\hfill\break}
\def\:{\,:\,}

\def\({\left(}                   
\def\){\right)}

\def\[{\left[}                   
\def\]{\right]}

\def\lan{\langle}
\def\ran{\rangle}

\def\ci{\subset}

\def\fy{\infty}

\def\del{\partial}

\def\lam{\lambda}

\def\Om{\Omega}   

\cl{\bf A compactness theorem of $n$-harmonic maps}
\bs
\cl{\au Changyou Wang}
\cl{Department of Mathematics, University of Kentucky}
\cl{Lexington, KY 40506, USA}

\bs
\nind{\bf Abstract}. {\it For $n\ge 3$, let $\Omega\ci \R^n$
be a bounded smooth domain and $N\ci \R^L$ be a compact smooth Riemannian submanifold
without boundary. Suppose that $\{u_n\}\ci W^{1,n}(\Omega, N)$ 
are weak solutions to the perturbed $n$-harmonic map equation (1.2), satisfying (1.3),
and $u_k\to u$ weakly in $W^{1,n}(\Omega, N)$. Then $u$ is an $n$-harmonic map. 
In particular, the space of $n$-harmonic maps 
is sequentially compact for the weak-$W^{1,n}$ topology.}

\bs
\nind{\S1} {\bf Introduction}
\ss
For $n\ge 2$, let $\Omega\ci \R^n$ be a bounded smooth domain,
and $N\ci \R^L$ be a compact smooth Riemannian manifold without boundary, 
isometrically embedded into the euclidean space $\R^L$ for some $L\ge 1$. 
For $2\le p\le n$, the Sobolev space $W^{1,p}(\Omega,N)$ is defined by
$$W^{1,p}(\Omega, N)=\{u=(u^1,\cdots, u^L)\in W^{1,p}(M,\R^L)|\ u(x)\in N 
\hbox{ for a.e. } x\in \Omega\}.$$
The Dirichlet $p$-energy functional $E_p: W^{1,p}(\Omega,N)\to \R$ is defined by
$$E_p(u)=\int_\Omega |\nabla u|^p\,dx= 
\int_\Omega (\sum_{\alpha=1}^n\lan {\del u\over \del x_\alpha},
{\del u\over\del x_\alpha}\ran)^{p\over 2}\,dx$$
where $\lan\cdot,\cdot\ran$ is the inner product of $\R^L$.

Recall that a map $u\in W^{1,p}(\Omega, N)$ is a $p$-harmonic map, 
if $u$ is a critical point of
$E_p(\cdot)$ on the space $W^{1,p}(\Omega, N)$, i.e.
$u$ satisfies the $p$-harmonic map equation:
$$-\hbox{div}(|\nabla u|^{p-2}\nabla u)=|\nabla u|^{p-2}A(u)(\nabla u,\nabla u) \eqno(1.1)$$
in the sense of distributions, 
where $\hbox{div}$ is the divergence operator on $\R^n$ and
$A(\cdot)(\cdot,\cdot)$ is the second fundamental form of $N\ci \R^L$.

Since the $p$-harmonic map equation (1.1) is an (degenerately) elliptic system
with critical nonlinearity in the first order derivatives, the analysis of
both the regularity problem for $p$-harmonic maps and the limiting behavior 
of weakly convergent sequences of $p$-harmonic maps are very interesting and extremely 
challenging.

In this paper, we are mainly interested in the compactness problem in
the weak topology of the space of $p$-harmonic maps. More precisely,
we are motivated by the following problem.
\ss
\nind{\bf Question A}. {\it For $n\ge 3$ and $2\le p\le n$, is any weak limit
$u$ in $W^{1,p}(\Omega,N)$ of a sequence of $p$-harmonic maps $\{u_k\}\ci
W^{1,p}(\Omega,N)$ is a $p$-harmonic map?} 

For $p=n=2$, the answer to question A is affirmative,
which follows from H\'elein's celebrated regularity theorem [H1]:
{\it any $2$-harmonic map from a Riemannian surface into any 
compact Riemannian manifold is smooth}.  

For $n\ge 3$, the answer to question A remains open in general cases,
although many people have made efforts to understand it. 

We mention some earlier results
in the direction. Schoen-Uhlenbeck [SU] ($p=2$),
Hardt-Lin [HL] and Luckhaus [L] ($p\not=2$) have shown that
{\it any weak limit $u\in W^{1,p}$ of a sequence of minimizing $p$-harmonic maps
is a strong limit and a minimizing $p$-harmonic map}. In particular, question A
is true for minimizing $p$-harmonic maps.

Without the minimality assumption, it is known that question A holds for target
manifolds $N$ with symmetry, such as $N=S^{L-1}$ is the unit
sphere in $\R^L$ (cf. Chen [C], Shatah [S], Evans [E2] \S5, and H\'elein [H2] \S2)
or $N={{\bf G}/{\bf H}}$ is a compact Riemannian homogeneous manifold 
(cf. Toro-Wang [TW]). Here the symmetry guarantees
the existence of Killing tangent vector fields on $N$, under which the nonlinearity
of the $p$-harmonic map equation (1.1) is the inner product of a gradient
and a divergence free vector field and hence belongs to the Hardy space ${\Cal H}^1(\R^n)$,
an improved subspace of $L^1(\R^n)$. 

For target manifolds $N$ without symmetry, 
the idea of the use of Coulomb moving frames, originally due to H\'elein [H1] and beautifully
explained in his book [H2], has played extremely important roles on the study of regularity 
of stationary $2$-harmonic maps into general target manifolds, through the works by 
H\'elein [H1] ($n=2$) and Bethuel [B2] ($n\ge 3$) (see also Evans [E1]).   
Roughly speaking, one can make suitable rotations from a smooth moving frames along $u^*TN$ 
to obtain a harmonic moving frame $\{e_\alpha\}$
(i.e. a minimizer of $\int |\lan de_\alpha, e_\beta\ran|^2$). It turns out 
that the nonlinearity of $2$-harmonic map equation via harmonic moving frames
enjoys the Jacobian structure partially. Although this is sufficient for the regularity
(and hence convergence) of stationary $2$-harmonic maps, it is not good enough for compactness
of weakly $2$-harmonic maps. On the other hand, in the study on existence of
wave maps in $\R^{2+1}$, Freire-M\"uller-Struwe [FMS1,2] have discovered that for the
class of wave maps enjoying the energy monotonicity inequality (e.g.
smooth wave maps) in $\R^{2+1}$, the concentration compactness method of Lions [L1,2],
in combination with the idea of Coulomb moving frames for wave maps and
some end-point analytic estimates, can be used to establish the compactness 
of wave maps in the class.   

When considering $p$-harmonic maps into general target manifolds $N$ for $p\not=2$, 
one may encounter the difficulty that is what might be the appropriate construction 
of Coulomb moving frames (e.g. neither minimizers of $\int |\lan de_\alpha,
e_\beta\ran|^p$ seem to fit the eqn. (1.1) well 
nor minimizers of $\int |\nabla u|^{p-2}|\lan de_\alpha, e_\beta\ran|^2$ seem to have
the $L^p$-estimate). However, we observe that, for $p=n$ case, Uhlenbeck's
construction of Coulomb gauges for Yang-Mills fields [U] can be adopted 
to obtain Coulomb moving frames along $u^*TN$ under the smallness of $E_n(u)$, 
see \S2 below for the detail and also Wang [W2,3] for its applications to biharmonic maps. 
With such a Coulomb moving frame along $u^*TN$, we are able to modify the analytic techniques
by [FMS2] to show the compactness of a Palais-Smale sequence (e.g. a
sequence of weakly convergent $n$-harmonic maps) of the Dirichlet $n$-energy functional $E_n$ on 
$W^{1,n}(\Omega, N)$.

In order to state our results, we first recall
\ss
\nind{\bf Definition}. {\it A sequence of maps $\{u_k\}\ci W^{1,n}(\Om,N)$ is 
a Palais-Smale sequence for the Dirichlet $n$-energy functional $E_n$
on $W^{1,n}(\Om,N)$, if the following two conditions hold:
\nind{(a)} $u_k\rightarrow u$ weakly in $W^{1,n}(\Om,N)$,
and {(b)} $E_n'(u_k)\rightarrow 0$ in $(W^{1,n}(\Om,N))^*$. Here $(W^{1,n}(\Om,N))^*$
is the dual of $W^{1,n}(\Om,N)$. }

Note that (b) is equivalent to that $u_k$ satisfies the
perturbed $n$-harmonic map equation
$$-\hbox{div}(|\nabla u_k|^{n-2}\nabla u_k)
=|\nabla u_k|^{n-2}A(u_k)(\nabla u_k, \nabla u_k)+\Phi_k,
\eqno(1.2)$$
in the sense of distributions, and satisfies 
$$\lim_{k\rightarrow\fy}\|\Phi_k\|_{(W^{1,n}(\Om,N))^*}=0. \eqno(1.3)$$

The question is whether any weak limit $u$ of a Palais-Smale sequence is an $n$-harmonic map. 
This is highly nontrivial, since $E_n$ is conformally invariant, i.e.
$E_n(u)=E_n(u\circ \Psi)$ for any $C^1$-conformal transformation $\Psi: \Om\to \Om$, and
the conformal group is non-compact and hence $E_n$ doesn't satisfy the Palais-Smale condition
(cf. [SaU]).  Our main result is
\ss
\nind{\bf Theorem B}. {\it For $n\ge 3$, assume that $\{u_k\}\ci W^{1,n}(\Om,N)$ 
satisfy the equation (1.2), (1.3), and converge weakly  to $u$ in $W^{1,n}(\Om,N)$, 
then $u\in W^{1,n}(\Om,N)$ is an $n$-harmonic map.}

We would like to remark that for $n=2$, theorem B has first been proven by Bethuel [B1],
later by Freire-M\"uller-Struwe [FMS2], and also by Wang [W1] with a method different
from both [B1] and [FMS2]. For $n\ge 3$, Hungerbhler [H] has obtained the existence of
global weak solutions to the $n$-harmonic map flow. Theorem B is applicable to the 
$n$-harmonic map flow by [H] at infinity time.

As a corollary, we confirm that question A is true for $p=n\ge 3$, i.e.
\ss
\nind{\bf Corollary C}. {\it For $n\ge 3$, assume that $\{u_k\}\ci W^{1,n}(\Om,N)$ 
are a sequence of $n$-harmonic maps converging
weakly to $u$ in $W^{1,n}(\Om,N)$, then $u$ is an $n$-harmonic map.}

The paper is written as follows. In \S2, we outline the construction of Coulomb
moving frames. In \S3, we first recall ${\Cal H}^1(\R^n)$-estimate for functions 
with Jacobian structure
by [CLMS], the duality between ${\Cal H}^1(\R^n)$ and $\hbox{BMO}(\R^n)$ by
[FS], and then give a proof of theorem B.

In this paper, we will use the following notations. For a ball $B=B_r(x)\ci\R^n$, denote
$\alpha B=B_{\alpha r}(x)$ for any $\alpha>0$. For $1\le i\le n$, denote
$\wedge^i(\R^n)$ as the i$^{\hbox{th}}$ wedge product of $\R^n$,  
$C^\fy(\R^n, \wedge^i(\R^n)$ as the space of smooth i$^{\hbox{th}}$ forms on $\R^n$,
and $W^{m,p}(\R^n, \wedge^i(\R^n)$ as the space of i$^{\hbox{th}}$ forms on $\R^n$
with $W^{m,p}(\R^n)$ coefficients, for nonnegative integers $m$ and $1<p<\fy$.
Denote by ${\Cal D}'(\Om)$ the dual of $C^\fy_0(\Om)$.
Denote $d$ as the exterior differentation operator on 
$\R^n$ and $\delta$ as the adjoint operator of $d$.

\bs
\nind{\S2} {\bf The construction of Coulomb moving frames}
\ss
This section is devoted to the construction of Coulomb moving frames 
along $u^*TN$, under the smallness condition on $E_n(u)$.

First recall that for any open set $U\ci \R^n$ and $u\in W^{1,n}(U,N)$, 
denote $u^*TN$ as the pull-back bundle of $TN$ by $u$ over $U$. 
Denote $l$ as the dimension of $N$,
we call $\{e_\alpha\}_{\alpha=1}^l$ a moving frame along $u^*TN$, 
if $\{e_\alpha(x)\}_{\alpha=1}^l$ 
forms an orthonormal base of $T_{u(x)}N$, the tangent space of $N$ at the point $u(x)$,
for a.e. $x\in U$. 

Now we have the perturbed $n$-harmonic map equation via a moving frame.
\ss
\nind{\bf Lemma 2.1}. {\it For $n\ge 3$, let $u\in W^{1,n}(\Omega,N)$ be a weak solution 
to the perturbed
$n$-harmonic map equation:
$$-\hbox{div}(|\nabla u|^{n-2}\nabla u)=|\nabla u|^{n-2}A(u)(\nabla u,\nabla u)+\Phi, 
\ \Phi\in (W^{1,n}(\Omega,N))^*. \eqno(2.1)$$
Suppose that $\{e_\alpha\}_{\alpha=1}^l$ is a moving frame along $u^*TN$. 
Then we have, for $1\le \alpha\le l$,
$$-\hbox{div}(\lan |\nabla u|^{n-2}\nabla u,e_\alpha\ran)
=\sum_{\beta=1}^l\lan |\nabla u|^{n-2}\nabla u, e_\beta\ran\lan\nabla e_\alpha, e_\beta\ran
+\lan \Phi, e_\alpha\ran \eqno(2.2)$$
in the sense of distributions.}
\ss
\nind{\bf Proof}. Note that for any $1\le \alpha\le l$ and  a.e. $x\in\Omega$, we have
$$\lan e_\alpha(x), A(u(x))(\nabla u(x),\nabla u(x))\ran =0$$
for
$e_\alpha(x)\in T_{u(x)}N$ and $A(u(x))(\nabla u(x), \nabla u(x))\perp T_{u(x)}N$.
Then straightforward calculations deduce (2.2) from (2.1).  \qed 

Now we establish a Coulomb moving frame along $u^*TN$, with the desired estimates
on its connection form. The constrction is inspired by an earlier result 
by the author in the context of biharmonic maps (cf. Wang [W2,3]) and 
Uhlenbeck's Coulomb gauge construction for Yang-Mills fields [U].

\ss
\nind{\bf Proposition 2.2}. {\it For $n\ge 3$ and any ball $B\ci\R^n$,
there exists an $\epsilon_0>0$ such that if $u\in W^{1,n}(2B,N)$
satisfies
$$\|\nabla u\|_{L^n(2B)}\le \epsilon_0 \eqno(2.3)$$
then there exists a Coulomb moving frame $\{e_\alpha\}_{\alpha=1}^l$
along $u^*TN$ in $W^{1,n}(B, \R^L)$ such that
its connection form $A=(\lan de_\alpha, e_\beta\ran)$ satisfies
$$\delta A=0 \ \hbox{ in } B; \ \  x\cdot A=0 \ \hbox{ on }\del B \eqno(2.4)$$
and
$$\|A\|_{L^n(B)}+\|\nabla A\|_{L^{n\over 2}(B)}
\le C\|\nabla u\|_{L^n(B)}^2. \eqno(2.5)$$}
\ss
\nind{\bf Proof}. Let us first assume  $u\in C^\fy(B,N)$.
Then $u^*TN|_B$ is a smooth vector bundle over the contractible
manifold $B$. Hence $u^*TN|_B$ is trivial 
and there exists a smooth moving
frame $\{{\bar e}_\alpha\}_{\alpha=1}^l$ along $u^*TN$
on $B$. Let ${\bf G}$ denote the gauge transformation group 
of $u^*TN$ consisting of maps $R:B\to {\bf SO}(l)$.
For any $R\in {\bf G}\cap W^{1,n}(B,{\bf SO}(l))$,
we rotate $\{{\bar e}_\alpha\}_{\alpha=1}^l$ to get another moving
frame $\{e_\alpha=\sum_{\beta=1}^l R^{\alpha\beta}{\bar e}_\beta\}_{\alpha=1}^l$.
Then we have
$$de_\alpha=\sum_{\beta=1}^l A^{\alpha \beta} e_\beta,  \ 1\le\alpha\le l,$$
where $(A^{\alpha\beta})=(\lan de_\alpha, e_\beta\ran)$ is the (matrix-valued) 
connection form of $u^*TN$.

Let $D$ denote the pull-back covariant derivatives on $u^*TN$. Then
the curvature equation of $u^*TN$ is given by: for $1\le p, q\le n$,
$$\eqalignno{D_pD_q e_\alpha-D_qD_pe_\alpha&=\sum_{\beta=1}^l(D_p(A^{\alpha \beta}_q
e_\beta)-
D_q(A^{\alpha \beta}_p e_\beta))\cr
&=\sum_{\beta=1}^l\{{\del}_p(A^{\alpha\beta}_q)-{\del}_q(A^{\alpha \beta}_p)
+\sum_{\gamma=1}^l(A^{\alpha\gamma}_q A^{\gamma\beta}_p
-A^{\alpha\gamma}_pA^{\gamma\beta}_q)e_\beta\}\cr
&=\sum_{\beta=1}^lF_{pq}^{\alpha\beta} e_\beta &(2.6)\cr}$$
or for short,
$$\del_p A_q-\del_q A_p+[A_p,A_q]=F_{pq}={\Cal R}^N(\del_p u,\del_q
u)\eqno(2.7)$$
where $A_p=(A_p^{\alpha\beta}=(\lan D_{\del\over \del x_p} e_\alpha, e_\beta\ran)$, 
and ${\Cal R}^N$ is the
curvature tensor of $TN$.
Here we have used the formula
$$\eqalignno{&D_pD_qe_\alpha-D_qD_p e_\alpha\cr
&=\sum_{\beta,\delta=1}^l\lan{\del u\over\del x_p}, e_\beta\ran
\lan{\del u\over\del x_q}, e_\delta\ran
u^*(D^N_{e_\beta}D^N_{e_\delta}e_\alpha-D^N_{e_\delta}D^N_{e_\beta}e_\alpha
-D^N_{[e_\beta,e_\delta]}e_\alpha)\cr
&=\sum_{\beta,\delta=1}^l\lan{\del u\over\del x_p}, e_\beta\ran
\lan{\del u\over\del x_q}, e_\delta\ran
u^*({\Cal R}^N_{e_\beta, e_\delta}(e_\alpha))&(2.8)\cr}$$
where $D^N$ is the Levi-Civita connection on $TN$.

For any $R\in {\bf G}\cap W^{1,n}(B, {\bf SO}(l))$, we know that
the connection form ${\bar A}=(\lan d{\bar e}_\alpha, {\bar e}_\beta\ran)$
of $\{{\bar e}_\alpha\}_{\alpha=1}^l$ and
the connection form $A=(\lan de_\alpha, e_\beta\ran)$ of 
$\{e_\alpha =\sum_{\beta=1}^lR^{\alpha\beta}{\bar e}_\beta\}_{\alpha=1}^l$ is
related by
$$A=R^{-1}dR+R^{-1} {\bar A}R. \eqno(2.9)$$
We also have the curvature $|F({\bar A})|(x)=|F(A)|(x)$ for a.e.
$x\in B$. Therefore 
the $L^{n\over 2}$-norm of curvature
$$\int_{B} |F(A)|^{n\over 2}\,dx=\int_{B}|F({\bar A})|^{n\over 2}\,dx$$
is invariant under gauge transformations.
Moreover, (2.7) implies that, for a.e. $x\in B$,
$$|F({\bar A})|(x)=|F(A)|(x)\le C\|{\Cal R}^N\|_{L^\fy}|\nabla u|^2(x)
\le C_N |\nabla u|^2(x).\eqno(2.10)$$
This implies 
$$\int_{B}|F(A)|^{n\over 2}\,dx\le C\int_{B}|\nabla u|^n\,dx. \eqno(2.11)$$

Now we use the condition (2.3) to approximate
$u\in W^{1,n}(B,N)$ by $C^\fy(B,N)$ as follows. Let $\phi:\R^n\to \R$
be a nonnegative, smooth radial mollifying function such that support $(\phi)\ci B_1$
and $\int_{R^n}\phi=1$. For $0<\epsilon<1$, let $\phi^\epsilon(x)
=\epsilon^{-n}\phi({x\over \epsilon})$ for $x\in \R^n$, and define
$$u^\epsilon(x)=\int_{\R^n}\phi^\epsilon(x-y)u(y)\,dy
=\int_{\R^n}\phi(y)u(x-\epsilon y)\,dy,  \ \forall x\in B. $$
For any $\epsilon\in (0, {1\over 2})$ and $x\in B$,
applying the modified Poincar\'e inequality to
$u_{x,\epsilon}(y)\equiv u(x-\epsilon y): B\to \R^L$, we
have
$$\int_{B}|u^\epsilon(x)-u_{x,\epsilon}|^n\,dy \le
C\int_{B}|\nabla u_{x,\epsilon}|^n\,dy
=C\int_{B_\epsilon(x)}|\nabla u|^n\,dy \le C\epsilon_0^n. \eqno(2.12)$$
Therefore we have, for any $\epsilon\in (0, {1\over 2})$,
$$\max_{x\in B}\hbox{dist}(u^\epsilon(x),N)\le C\epsilon_0$$
so that $u^\epsilon(B)\ci N_{C\epsilon_0}$. Since
the nearest point projection map $\Pi:N_{C\epsilon_0}\to N$ is smooth for
sufficiently small $\epsilon_0>0$, we have
$u_\epsilon=\Pi\circ u^\epsilon\in C^\fy(B, N)$,
$u_\epsilon\in C^\fy(B,N)$, and 
$u_\epsilon\rightarrow u$ strongly in $W^{1,n}(B,N)$
as $\epsilon\to 0$, and
$$\sup_{0<\epsilon<{1\over 2}}\|\nabla u_\epsilon\|_{W^{1,n}(B)}
\le C\|\nabla u\|_{W^{1,n}(2B)}\le C\epsilon_0. \eqno(2.13)$$
For $\epsilon\in (0, {1\over 2})$, since 
$u_\epsilon^*TN|_B$ are trivial, 
there exist smooth moving frames
$\{{\bar e}_\alpha^\epsilon\}_{\alpha=1}^n$ along $u_\epsilon^*TN$
over $B$. Moreover (2.13) and (2.11) imply that
the curvature of the connections ${\bar A}_\epsilon
=(\lan D{\bar e}_\alpha^\epsilon, {\bar e}_\beta^\epsilon\ran)$ satisfies
$$\|F({\bar A}_\epsilon)\|_{L^{n\over 2}(B)}\le C\|\nabla u_\epsilon\|_{W^{1,n}(B)}^2
\le C\|\nabla u\|_{W^{1,n}(2B)}^2\le C\epsilon_0^2. \eqno(2.14)$$
Since $\{{\bar e}^\epsilon_\alpha\}_{\alpha=1}^l$ are smooth frames and their
connections ${\bar A}_\epsilon$ satisfy (2.14) with sufficiently small
$\epsilon_0>0$, we can apply Uhlenbeck [U]
to conclude that there are $R_\epsilon: B\to {\bf SO}(l)$
satisfying $\nabla R_\epsilon\in L^n(B)$ and $\nabla^2 R_\epsilon\in L^{n\over 2}(B)$
such that the connection form
$A_\epsilon\equiv R_\epsilon^{-1}dR_\epsilon+R^{-1}_\epsilon {\bar A}_\epsilon R_\epsilon
$
of the moving frame $\{e_\alpha^\epsilon(x)\equiv
\sum_{\beta=1}^l R_\epsilon^{\alpha \beta}{\bar e}^\epsilon_\beta\}_{\alpha=1}^l$
satisfy
$$\eqalignno{\delta A_\epsilon = 0 \ \hbox{ in } B, \  \ \ &
             x\cdot A_\epsilon =0, \hbox{ on }\del B, &(2.15)\cr
             \|A_\epsilon\|_{L^n(B)}+\|\nabla A_\epsilon\|_{L^{n\over 2}(B)}
&\le C\|F({\bar A}_\epsilon)\|_{L^{n\over 2}(B)}\le C\epsilon_0. &(2.16)\cr}$$

We now estimate $\|\nabla e_\alpha^\epsilon\|_{L^n(B)}$ 
for $1\le \alpha\le l$. 

For $y\in N$, let $P^\perp(y)=y-\nabla \Pi(y): \R^L\to (T_y N)^\perp$
denote the orthogonal projection from map $\R^L$ to the normal space $(T_yN)^\perp$.
Then we have
$$\nabla e_\alpha^\epsilon=\sum_{\beta=1}^l\lan \nabla e_\alpha^\epsilon, e_\beta^\epsilon\ran
e_\beta^\epsilon+P^\perp(u_\epsilon)(\nabla e_\alpha^\epsilon)
=\sum_{\beta=1}^l \lan \nabla e_\alpha^\epsilon, e_\beta^\epsilon\ran e_\beta^\epsilon
-A(u_\epsilon)(e_\alpha^\epsilon,\nabla u_\epsilon) \eqno(2.17)$$
where we have used
$$P^\perp(u_\epsilon)(\nabla e_\alpha^\epsilon)=-\nabla (P^\perp(u_\epsilon))
(e_\alpha^\epsilon)=-A(u_\epsilon)(e_\alpha^\epsilon,\nabla u_\epsilon)$$
for $P^\perp(u_\epsilon)(e_\alpha^\epsilon)=$. 
Therefore we have
$$|\nabla e_\alpha^\epsilon|(x)\le C(|A_\epsilon|+|\nabla u_\epsilon|)(x), \ \ \hbox{
for a.e. }x\in B. \eqno(2.18)$$
This, combined with (2.13) and (2.16), yields
$$\sum_{\alpha=1}^l\|\nabla e_\alpha^\epsilon\|_{L^n(B)}
\le C(\|A_\epsilon\|_{L^n(B)}+\|\nabla u_\epsilon\|_{L^n(B)})
\le C\|\nabla u\|_{L^n(B)}. \eqno(2.19)$$
Therefore, after passing to subsequences,
we can assume that
$e_\alpha^\epsilon\rightarrow e_\alpha$ weakly in $W^{1,n}(B)$, strongly
in $L^n(B)$, and a.e. in $B$. Since $u_\epsilon\rightarrow u$ strongly
in $W^{1,n}(B)$, we have that $\{e_\alpha\}_{\alpha=1}^l,\ci W^{1,n}(B)$,
is a moving frame along $u^*TN$. Moreover, (2.16) implies that 
$A_\epsilon\rightarrow A \equiv(\lan d e_\alpha, e_\beta\ran)$, the connection form of 
$\{e_\alpha\}_{\alpha=1}^l$,
weakly in $W^{1,{n\over 2}}(B)$.
Hence, by taking $\epsilon$ into zero, (2.15) and (2.16) imply that
$A$ satisfies (2.4) and (2.5).
The proof of Proposition 2.2 is complete.  \qed
\bs
\bs
\nind{\S3} {\bf Proof of theorem B}

This section is devoted to the proof of theorem B.  First we recall
some basic facts on the Hardy space ${\Cal H}^1(\R^n)$ and the BMO space
$\hbox{BMO}(\R^n)$. 

Recall that $f\in L^1(\R^n)$ belongs to the Hardy space ${\Cal H}^1(\R^n)$ if
$$f_*:=\sup_{\epsilon>0}|\phi_\epsilon*f|\in L^1(\R^n)$$
where $\phi_\epsilon(x):=\epsilon^{-n}\phi({x\over\epsilon})$ for
a fixed nonnegative $\phi\in C_0^\fy(\R^n)$ with $\int_{\R^n}\phi\,dy=1$.
Note that ${\Cal H}^1(\R^n)$ is a Banach space with the norm
$$\|f\|_{{\Cal H}^1(\R^n)}:=\|f\|_{L^1(\R^n)}+\|f_*\|_{L^1(\R^n)}.$$
An important property of $f\in {\Cal H}^1(\R^n)$ is the cancellation
identity $\int_{\R^n}f\,dy=0$ (cf. Fefferman-Stein [FS]).

Recall also that  $f\in L^1_{\hbox{loc}}(\R^n)$ belongs to the BMO space $\hbox{BMO}(\R^n)$
(cf. John-Nirenberg [JN]), if
$$\|f\|_{\hbox{BMO}(\R^n)}
:=\sup\{{1\over |B|}\int_{B}|f-f_B|\,dy: \hbox{ any ball }B\ci \R^n\}<\fy$$
where $f_B={1\over |B|}\int_B f\,dy$ is the average of $f$ over $B$. By
the Sobolev inequality we have $W^{1,n}(\R^n)\ci \hbox{BMO}(\R^n)$ and
$$\|f\|_{\hbox{BMO}(\R^n)}\le C\|\nabla f\|_{L^n(\R^n)}. \eqno(3.1)$$

The celebrated theorem of Fefferman-Stein [FS] says that the dual of ${\Cal H}^1(\R^n)$
is $\hbox{BMO}(\R^n)$. Moreover
$$|\int_{\R^n} fg\,dy|\le C\|f\|_{{\Cal H}^1(\R^n)}\|g\|_{\hbox{BMO}(\R^n)}. \eqno(3.2)$$

Now we recall an important result of 
Coifman-Lions-Meyer-Semmes [CLMS], see also [E1].
\ss
\nind{\bf Proposition 3.1} ([CLMS]). {\it For any $1<p<\fy$, denote $p'={p\over p-1}$.
Let $f\in W^{1,p}(\R^n)$, $g\in W^{1,p'}(\R^n, \wedge^1(\R^n))$, and
$h\in W^{1,n}(\R^n)$. Then
$df\cdot \delta g\in {\Cal H}^1(\R^n)$ and
$$\|df\cdot \delta g\|_{{\Cal H}^1(\R^n)}\le C\|\nabla f\|_{L^p(\R^n)}
\|\nabla g\|_{L^{p'}(\R^n)}. \eqno(3.3)$$
In particular, we have 
$$|\int_{\R^n} \lan df\cdot \delta g, h\ran \,dy| 
\le C\|\nabla f\|_{L^p(\R^n)}\|\nabla g\|_{L^{p'}(\R^n)}\|\nabla h\|_{L^n(\R^n)}. \eqno(3.4)$$}

We also recall the following pointwise convergence result, which is essentially
due to Hardt-Lin-Mou [HLM] (see also [F]).
\ss
\nind{\bf Lemma 3.2} ([HLM]). {\it Suppose that $\{u_k\}\ci W^{1,n}(\Omega, \R^L)$ 
are weak solutions to
$$-\hbox{div}(|\nabla u_k|^{n-2}\nabla u_k)=f_k+\Phi_k, \eqno(3.5)$$
where $f_k\to 0$ in $L^1(\Omega, \R^L)$, and $\Phi_k\to 0$ in 
$(W^{1,n}(\Omega, \R^L))^*$.
Assume that $u_k\to u$ weakly in $W^{1,n}(\Omega, \R^L)$. Then, after taking
possible subsequences, we have $\nabla u_k\rightarrow \nabla u$ a.e. in $\Omega$.
In particular, $\nabla u_k\to \nabla u$ strongly in $L^q(\Omega, \R^L)$
for any $1\le q<n$.}

After these preparations, we are ready to give a proof of theorem B. It turns
out the crucial step is to show the following weak compactness under the smallness
condition on $E_n$.  
\ss
\nind{\bf Lemma 3.3} ($\epsilon$-weak compactness). {\it For any $n\ge 3$, there exists
an $\epsilon_1>0$ such that if $\{u_k\}\ci W^{1,n}(2B,N)$ satisfy both the equation (1.2)
and the condition (1.3) with $\Omega$ replaced by $2B$,
$u_k\to u$ weakly in $W^{1,n}(2B,N)$, and satisfy
$$\int_{2B}|\nabla u_k|^n\,dx \le \epsilon_1^n, \ \ \forall k\ge 1. \eqno(3.6)$$
Then $u\in W^{1,n}(B,N)$ is an $n$-harmonic map.}
\ss
\nind{\bf Proof}. For the convenience, we will write both equation (1.1) and
(1.2) by using $d$ and $\delta$ from now on. 

Let $\epsilon_1>0$ be the same constant as in Proposition 2.2.
Then we have that for any $k\ge 1$ there is a Coulomb moving frame
$\{e_\alpha^k\}_{\alpha=1}^l$ along $u_k^*TN$ such that the connection form
$A_k=(\lan de_\alpha^k, e_\beta^k\ran)$ satisfies
$$\delta A_k=0\ \hbox{ in }B; \ \  x\cdot A_k =0 \hbox{ on }\del B \eqno(3.7)$$
and
$$\|A_k\|_{L^n(B)}+\|\nabla A_k\|_{L^{n\over 2}(B)}
\le C\|\nabla u_k\|_{L^n(B)}^2. \eqno(3.8)$$
Moreover, similar to (2.19),  we have
$$\max_{\alpha=1}^l\|\nabla e_\alpha^k\|_{L^n(B)}\le C\|\nabla u_k\|_{L^n(B)}\le C\epsilon_1,
\ \ \forall k\ge 1. \eqno(3.9)$$
Therefore we may assume, after passing to subsequences, that
$e_\alpha^k\rightarrow e_\alpha$ weakly in $W^{1,n}(B,\R^L)$ and strongly in $L^n(B,\R^L)$, 
$A_k\rightarrow A$ weakly in $W^{1,{n\over 2}}(B)$ and strongly in $L^{n\over 2}(B)$.
It is easy to see that $\{e_\alpha\}_{\alpha=1}^l$
is a moving frame along $u^*TN$, and  $A=(\lan de_\alpha,e_\beta\ran)$ 
satisfies
$$\delta A=0 \hbox{ in } B; \ \  \ x\cdot A=0 \hbox{ on }\del B, \eqno(3.10)$$
and
$$\|A\|_{L^n(B)}+\|\nabla A\|_{L^{n\over 2}(B)}\le C\liminf_{k}\|\nabla u_k\|_{L^n(B)}^2\le
C\epsilon_1^2. \eqno(3.11)$$

Using these moving frames, Lemma 2.1 yields that for any $1\le\alpha\le l$
$$-\delta(\lan |du_k|^{n-2}du_k, e_\alpha^k\ran)
=\sum_{\beta=1}^l \lan |du_k|^{n-2}du_k, e_\beta^k\ran \cdot\lan de_\alpha^k, e_\beta^k\ran
+\lan \Phi_k, e_\alpha^k\ran.\eqno(3.12)$$

It follows from Lemma 3.2 that we can assume that $\nabla u_k\rightarrow \nabla u$
strongly in $L^q(\Omega)$ for any $1\le q<n$. Therefore we have
$$|du_k|^{n-2}du_k\to |du|^{n-2} du, \  \hbox{  weakly in }
L^{n\over n-1}(2B). \eqno(3.13)$$
This implies 
$$-\delta(\lan |du_k|^{n-2}du_k, e_\alpha^k\ran)\rightarrow
-\delta(\lan |du|^{n-2}du, e_\alpha\ran), \hbox{ in } {\Cal D}'(B)\eqno(3.14)$$
as $k\rightarrow\fy$, for all $1\le \alpha\le l$.

It is readily seen that for any $\phi\in C^\fy_0(B)$ we have
$$|\lan \Phi_k, e_\alpha^k \phi\ran_{\{(W^{1,n})^*, W^{1,n}\}}|
\le \|\Phi_k\|_{(W^{1,n}(B,N))^*}\|e_\alpha^k\phi\|_{W^{1,n}(B)} \rightarrow 0,
\hbox{ as } k\rightarrow\fy. \eqno(3.15)$$
In order to prove that $u$ is an $n$-harmonic map, it
suffices to prove that for any $1\le\alpha,\beta\le l$
$$\lan |du_k|^{n-2}du_k, e_\beta^k\ran \cdot\lan de_\alpha^k, e_\beta^k\ran 
\rightarrow \lan |du|^{n-2}du, e_\beta\ran \lan de_\alpha, e_\beta\ran,
\hbox{ in } {\Cal D}'(B). \eqno(3.16)$$
 
To prove (3.16), we first let $\bar u_k\in W^{1,n}(\R^n, \R^L)$ and
$\overline{e_\alpha^k}\in W^{1,n}(\R^n,\R^L)$ be the extensions of $u_k$
and $e_\alpha^k$ from $B$ respectively such that
$$\|\nabla \bar u_k\|_{L^n(\R^n)}\le C\|\nabla u_k\|_{L^n(B)}, 
\ \|\nabla (\overline {e_\alpha^k})\|_{L^n(\R^n)}\le C\|\nabla e_\alpha^k\|_{L^n(B)}.
\eqno(3.17)$$
For $\lan|d\bar u_k|^{n-2}d\bar u_k, \overline{e_\beta^k}\ran\in L^{n\over n-1}(\R^n, \wedge^1(\R^n))$, 
the Hodge decomposition theorem (cf. Iwaniec-Martin [IM]) implies 
that there are $f_\beta^k\in W^{1, {n\over n-1}}(\R^n)$ and 
$g_\beta^k\in W^{1,{n\over n-1}}(\R^n, \wedge^2(R^n))$ such that $dg_\beta^k=0$,
$$
\lan |d\bar u_k|^{n-2}d\bar u_k, \overline{e_\beta^k}\ran=df_\beta^k +\delta  g_\beta^k, \eqno(3.18)$$
and
$$\|\nabla f_\beta^k\|_{L^{n\over n-1}(\R^n)}+\|\nabla g_\beta^k\|_{L^{n\over n-1}(\R^n)}
\le C\|\nabla u_k\|_{L^n(B)}^{n-1}. \eqno(3.19)$$
It follows from (3.19) that we may assume 
$f_\beta^k\rightarrow f_\beta, g_\beta^k\rightarrow g_\beta$ weakly in 
$W^{1,{n\over n-1}}_{\hbox{loc}}(\R^n)$. Therefore, by taking $k$ to infinity, 
(3.18) implies
$$\lan |du|^{n-2}du, e_\beta\ran =df_\beta+\delta g_\beta;\   \ dg_\beta=0, \ \hbox{ in } B. \eqno(3.20)$$
Moreover, (3.18) gives
$$\lan |du_k|^{n-2}du_k, e_\beta^k\ran\cdot\lan de_\alpha^k, e_\beta^k\ran
=df_\beta^k \cdot\lan de_\alpha^k, e_\beta^k\ran 
+\delta g_\beta^k\cdot\lan de_\alpha^k, e_\beta^k\ran, \hbox{ in } B .\eqno(3.21)$$
Since $df_\beta^k\rightarrow df_\beta$ weakly in $L^{n\over n-1}(B)$,
$\lan de_\alpha^k, e_\beta^k\ran \rightarrow \lan de_\alpha, e_\beta\ran$ weakly
in $L^n(B)$, and $\delta\lan de_\alpha^k, e_\beta^k\ran =0$ in $B$, 
we can apply the Div-Curl Lemma (cf. Evans [E2] page 53) to conclude 
$$df_\beta^k \cdot\lan de_\alpha^k, e_\beta^k\ran 
\rightarrow df_\beta \cdot \lan de_\alpha, e_\beta\ran,
\hbox{ in }{\Cal D}'(B). \eqno(3.22)$$
In fact, (3.22) follows directly from the integrations by parts:
for any $\phi\in C^\fy_0(B)$, 
$$\eqalignno{&\int_{\R^n} df_\beta^k\cdot\lan de_\alpha^k, e_\beta^k\ran \phi\,dx
=-\int_{\R^n}f_\beta^k\lan de_\alpha^k, e_\beta^k\ran\cdot d\phi\,dx\cr
&\rightarrow -\int_{\R^n}f_\beta\lan de_\alpha, e_\beta\ran\cdot d\phi\,dx
=\int_{\R^n}df_\beta\cdot \lan de_\alpha, e_\beta\ran \phi\cr}$$
as $k\to\fy$. Here
we have used both (3.7) and (3.10), i.e.
$\delta\lan de_\alpha^k, e_\beta^k\ran=\delta\lan de_\alpha, e_\beta\ran=0, \hbox{ in } B$.

Now we need the compensated compactness result
(cf. Lions [L1,2]), which was developped by Freire-M\"uller-Struwe [FMS1,2]
in the context of wave maps on $\R^{2+1}$.
\ss
\nind{\bf Lemma 3.4}. {\it Under the same notations. After taking possible
subsequences, we have
$$\delta g_\beta^k\cdot \lan de_\alpha^k, e_\beta^k\ran
\rightarrow  \delta g_\beta\cdot
\lan de_\alpha, e_\beta\ran +\nu, \hbox{ in }B \eqno(3.23)$$
where $\nu$ is a signed Radon measure given by
$$\nu=\sum_{j\in J}a_j\delta_{x_j}\eqno(3.24)$$
where $J$ is at most countable, $a_j\in \R$, $x_j\in B$, and
$\sum_{j\in J}|a_j|<+\fy.$} 
\ss
\nind{\bf Proof}. For the simplicity,  we only outline
a proof based on suitable modifications of [FMS2].
 
First we observe that 
$$\eqalignno{&\delta g_\beta^k\cdot\lan de_\alpha^k, e_\beta^k\ran
-\delta g_\beta\cdot\lan de_\alpha, e_\beta\ran\cr
&=\delta(g_\beta^k-g_\beta)\cdot\lan d(e_\alpha^k-e_\alpha), e_\beta^k\ran
+\delta g_\beta\cdot\lan d(e_\alpha^k-e_\alpha), e_\beta^k\ran\cr
&+(\delta g_\beta^k\cdot\lan de_\alpha, e_\beta^k\ran 
-\delta g_\beta\cdot\lan de_\alpha, e_\beta\ran)\cr
&=\delta (g_\beta^k-g_\beta)\cdot\lan d(e_\alpha^k-e_\alpha), e_\beta^k\ran+I_k+II_k.\cr}$$
The dominated convergence theorem implies
$$I_k, II_k \rightarrow 0, \hbox{ in }L^1(B),  \ \hbox{ as } k\rightarrow\fy.$$ 
Therefore (3.23) and (3.24) is equivalent to
$$\delta(g_\beta^k-g_\beta)\cdot
\lan d(e_\alpha^k-e_\alpha), e_\beta^k\ran
\rightarrow \nu \eqno(3.25)$$
where $\nu$ is the Radon measure given by (3.24).

Since $|\nabla (e_\alpha^k-e_\alpha)|^n$,
$|\nabla (g_\beta^k-g_\beta)|^{n\over n-1}$ are
bounded in $L^1(B)$, we may assume, after taking subsequences, that
there is a nonnegative Radon measure $\mu$ on $B$ such that
$$(\sum_{\alpha=1}^l|\nabla(e_\alpha^k-e_\alpha)|^n
+\sum_{\beta=1}^l|\nabla (g_\beta^k-g_\beta)|^{n\over 
n-1})\,dx
\rightarrow \mu$$
as convergence of Radon measures on $B$. 

Let ${\Cal S}=\{x\in B: \mu(\{x\})\equiv \lim_{r\rightarrow 0}\mu(B_r(x))>0\}.$
Then it follows from $\mu(B)<+\fy$ that ${\Cal S}$ is at most a countable
set. Now we want to show 
$$\hbox{supp}(\nu)\ci {\Cal S}. \eqno(3.26)$$
It is easy to see that (3.26) yields (3.24) and hence the conclusion of Lemma 3.4.

To see (3.26), we proceed as follows.
For $\phi\in C^\fy_0(B)$, we have
$$\eqalignno{&\lan\nu, \phi^3\ran 
=\lim_{k\to\fy}\int_{\R^n}\phi\delta(g_\beta^k-g_\beta)\cdot
\lan \phi d({e_\alpha^k-e_\alpha}), \phi e_\beta^k\ran \,dx\cr
&=\lim_{k\to\fy}\int_{\R^n}[\delta(\phi(g_\beta^k-g_\beta))-d\phi\cdot (g_\beta^k-g_\beta)]
\cdot\lan [d(\phi(e_\alpha^k-e_\alpha))-(e_\alpha^k-e_\alpha)d\phi], \phi e_\beta^k\ran\,dx\cr
&=\lim_{k\rightarrow\fy}\int_{\R^n}
\delta(\phi(g_\beta^k-g_\beta))\cdot\lan d(\phi(e_\alpha^k-e_\alpha)), \phi e_\beta^k\ran\,dx &
(3.27)\cr}$$
where we have used 
$$\lim_{k\to\fy}\int_{\R^n}[(g_\beta^k-g_\beta)d\phi\cdot\lan \phi d(e_\alpha^k-e_\alpha), \phi
e_\beta^k\ran
-\delta(\phi(g_\beta^k-g_\beta))\cdot\lan (e_\alpha^k-e_\alpha)d\phi, \phi e_\beta^k\ran]\,dx=0.$$

Note that Proposition 3.1 implies
$H_k\equiv \delta(\phi(g_\beta^k-g_\beta))\cdot d(\phi(e_\alpha^k-e_\alpha))$ is bounded
in ${\Cal H}^1(\R^n)$,
and (3.22) implies $H_k\rightarrow 0$ in ${\Cal D}'(\R^n)$. Therefore we have that
$H_k\rightarrow 0$ weak$^*$ in ${\Cal H}^1(\R^n)$. On the other hand, since
$\phi e_\beta\in W^{1,n}(\R^n)$, we have $\phi e_\beta\in \hbox{ VMO}(\R^n)$, where
$\hbox{VMO}(\R^n)\ci \hbox{BMO}(\R^n)$ is the closure of $C^\fy_0(\R^n)$ in the $\hbox{BMO}$ norm.
It is well-known [FS] that the dual of $\hbox{VMO}(\R^n)$ is ${\Cal H}^1(\R^n)$. Hence we have
$$\lim_{k\to\fy}\int_{\R^n} 
\delta(\phi(g_\beta^k-g_\beta))\cdot \lan d(\phi(e_\alpha^k-e_\alpha)), \phi e_\beta\ran\,dx
=0. \eqno(3.28)$$
Putting (3.28) together with (3.27) and applying (3.4), we have
$$\eqalignno{&|\lan\nu, \phi^3\ran|\cr
&\le C\lim_{k\rightarrow\fy}\|\nabla(\phi (e_\beta^k-e_\beta))\|_{L^n(\R^n)}
\|\nabla(\phi ({e_\alpha^k-e_\alpha}))\|_{L^n(\R^n)}
\|\nabla(\phi({g_\beta^k-g_\beta}))\|_{L^{n\over n-1}(\R^n)}\cr
&\le C\lim_{k\rightarrow\fy}\{[\|\phi\nabla (e_\beta^k-e_\beta)\|_{L^n(\R^n)}+\|\nabla \phi\|_{L^\fy}
\|e_\beta^k-e_\beta\|_{L^n(B)}]\cr
&\ \ \cdot [\|\phi\nabla(e_\alpha^k-e_\alpha)\|_{L^n(\R^n)}+\|\nabla \phi\|_{L^\fy}
\|e_\alpha^k-e_\alpha\|_{L^n(B)}]\cr
&\ \ \cdot[\|\phi\nabla(g_\beta^k-g_\beta)\|_{L^{n\over n-1}(\R^n)}+\|\nabla\phi\|_{L^\fy}
\|g_\beta^k-g_\beta\|_{L^{n\over n-1}(B)}]\}\cr
&\le C(\lan \mu, \phi^n\ran)^{1\over n}(\lan \mu, \phi^n\ran)^{1\over n}(\lan \mu,
\phi^{n\over n-1}\ran)^{1\over n-1}&(3.29)\cr}$$
where we have used
$$\lim_{k\to\fy}(\|e_\alpha^k-e_\alpha\|_{L^n(B)}+\|g_\beta^k-g_\beta\|_{L^{n\over n-1}(B)})=0.$$

By choosing $\phi_i\in C^\fy_0(B)$ such that $\phi_i\to \lam_{B_r(y)}$,
the characteristic function of a ball $B_r(y)$, we then have
$$\nu(B_r(y))\le C\mu(B_r(y))^{n+1\over n}. \eqno(3.30)$$
Therefore $\nu$ is absolutely continuous with respect to
$\mu$. Moreover,  for any $y\notin {\Cal S}$, the 
Radon-Nikodyn derivative
$${d\nu\over d\mu}(y)=\lim_{r\to 0}{\nu(B_r(y))\over \mu(B_r(y))}
\le C\lim_{r\to 0}\mu(B_r(y))^{1\over n}=0.$$
Therefore the support of $\nu$ is contained in ${\Cal S}$. This
proves (3.26) and hence (3.24).  The proof of Lemma 3.4 is  complete. \qed

Now we return to the proof of Lemma 3.3. By putting (3.14), (3.20), (3.22), and (3.23) together,
we have, for any $1\le\alpha\le l$,
$$-\delta(\lan |du|^{n-2}du,e_\alpha\ran)=\sum_{\alpha=1}^l\lan |du|^{n-2}du, e_\beta\ran \cdot \lan
de_\alpha, e_\beta\ran+\sum_{j\in J}a_j \delta_{x_j} \eqno(3.31)$$
where $J$ is at most countable, $a_j\in\R$, $x_j\in B$, and $\sum_{j\in J}|a_j|<+\fy$.

In order to conclude that $u$ is an $n$-harmonic map, one has to show that $a_j=0$
for all $j\in J$. This can be done by the following simple capacity argument.

Fix any $j\in J$, let $\phi\in C_0^\fy(\R^n)$ be such that $\phi\equiv 1 $ on ${1\over 2}B$ and
$\phi\equiv 0$ off $B$. For $k\ge 1$, let $\phi_k(x)=\phi(k(x-x_l))$.
Then, as $k\to\fy$,
$$\phi_k\rightarrow 0 \hbox{ weakly in } W^{1,n}(B), \phi_k\rightarrow 0 \hbox{ boundedly a.e.}$$
Now test (3.31) with $\phi_k$, let $k\to\fy$ and use the dominated convergence theorem
to deduce $a_j=0$. Therefore  (3.31) implies
$$-\delta(\lan |du|^{n-2}du,e_\alpha\ran)=\sum_{\alpha=1}^l\lan |du|^{n-2}du, e_\beta\ran \cdot \lan
de_\alpha, e_\beta\ran. \eqno(3.32)$$
It is easy to see that (3.32) is equivalent to (1.1), for $\{e_\alpha\}_{\alpha=1}^l$
is a moving frame along $u^*TN$. The proof of Lemma 3.3 is complete.   \qed

\ms
Based on Lemma 3.3, we can give a proof of theorem B as follows.
\ms
\nind{\bf Proof of theorem B}.

Since $|\nabla u_k|^n$ is bounded in $L^1(\Omega)$, we may assume, after passing to subsequences,
that there is a nonnegative Radon measure $\mu$ on $\Omega$ such that
$$|\nabla u_k|^n\,dx\rightarrow \mu$$
as convergence of Radon measures. Let $\epsilon_1>0$ be the same constant as in Lemma 3.3
and define $\Sigma\ci\Omega$ by
$$\Sigma=\{x\in\Omega: \mu(\{x\})\ge \epsilon_1^n\}.$$
Then $\Sigma$ is a finite subset and
$$|\Sigma|\le C\epsilon_1^{-n}, \ \ C\equiv\limsup_{k\to\fy}
\int_\Om |\nabla u_k|^n\,dx<+\fy.$$

For any $x_0\in\Omega\setminus \Sigma$, 
there exists an $r_0>0$ such that $\mu(B_{4r_0}(x_0))<\epsilon_1^n$. Since
$$\limsup_{k\to\fy}\int_{B_{2r_0}(x_0)}|\nabla u_k|^n\,dx\le \mu(B_{4r_0}(x_0)),$$
we can assume that there exists $k_0\ge 1$ such that
$\int_{B_{2r_0}(x_0)}|\nabla u_k|^2\,dx\le \epsilon_1^n, \ \ \forall k\ge k_0.$
Therefore Lemma 3.3 implies that $u$ 
is an $n$-harmonic map in $B_{r_0}(x_0)$. Since $x_0\in\Om\setminus\Sigma$ is
arbitrary, we conclude that $u$ is an $n$-harmonic map in $\Om\setminus\Sigma$. 

To show $u$ is an $n$-harmonic map in $\Om$, observe that $\hbox{Cap}_n(\Sigma)=0$
(cf. [EG]). Therefore there are
a sequence $\{\eta_i\}$ of functions on $\R^n$ such that $0\le \eta_i\le 1$,
$$\Sigma\ci\hbox{int}\{\eta_i=1\},\ \
\int_{\R^n}|\nabla\eta_i|^n\,dx\to 0,\ \ \eta_i\rightarrow 0 \hbox{ boundedly a.e.} \eqno(3.33)$$
Then, for any $\phi\in C^\fy_0(\Om,\R^L)$, we have
$$\int_{\Om}|du|^{n-2}du\cdot d\phi\,dx
=\int_\Om |du|^{n-2}du\cdot d((1-\eta_i)\phi)\,dx 
+\int_{\Om}|du|^{n-2}du\cdot (\eta_i d\phi+d\eta_i \phi)\,dx.$$
From (3.33) we  conclude that 
$(1-\eta_i)\phi$ is compactly supported in $\Om\setminus\Sigma$ and therefore,
$$\int_\Om |du|^{n-2}du\cdot d((1-\eta_i)\phi)\,dx 
=\int_\Om |du|^{n-2}A(u)(du, du)(1-\eta_i)\phi\,dx.$$
By the dominated convergence theorem and (3.33), we have
$$\lim_{i\to\fy}\int_\Om |du|^{n-2}A(u)(du, du)\eta_i\phi\,dx=0, \ 
\ \lim_{i\to\fy}\int_{\Om}|du|^{n-2}du\cdot \eta_i d\phi\,dx=0.$$
Applying (3.34) we have
$$|\int_\Om |du|^{n-2}du \cdot d\eta_i \phi\,dx|
\le \|\phi\|_{L^\fy(\Om)}(\int_\Om |Du|^n\,dx)^{n-1\over n}
(\int_{\R^n}|D\eta_i|^n\,dx)^{1\over n}\rightarrow 0$$
as $i\to\fy$. Therefore we have
$$\int_{\Om}|du|^{n-2}du\cdot d\phi\,dx
=\int_\Om |du|^{n-2}A(u)(du, du)\phi\,dx, \ \forall \phi\in C^\fy_0(\Om,\R^L).$$
\qed

\newpage
\cl{\bf REFERENCES}
\bs
\nind{[B1]}  F. Bethuel,
{\it Weak limits of Palais-Smale sequences for a class of critical functionals.}
Calc. Var. \& PDE 1 (1993), no. 3, 267--310.

\ss      
\nind{[B2]}  F. Bethuel, {\it On the singular set of stationary harmonic maps}. 
Manuscripta Math. 78 (1993),417-443 (1993).

\ss
\nind{[C]}  Y. M. Chen, {\it The weak solutions to the evolution problems of harmonic maps}. 
Math. Z. 201 (1989), no. 1, 69--74. 

\ss
\nind{[CLMS]}
R. Coifman, P. Lions, Y. Meyer, S. Semmes, {\it Compensated compactness and Hardy spaces}. 
J. Math. Pures Appl. 72 (1993), 247--286. 

\ss
\nind{[E1]}  L. C. Evans, {\it Partial regularity for stationary harmonic maps into spheres.} 
Arch. Rational Mech. Anal. 116 (1991), 101-113.

\ss
\nind{[E2]}  L. C. Evans, {Weak convergence methods for nonlinear partial differential equations.} 
CBMS Regional Conference Series in Mathematics, 74, 1990. 

\ss
\nind{[EG]}  L. C. Evans, R. Gariepy, {Measure theory and fine properties of functions}. 
Studies in Advanced Mathematics. CRC Press, Boca Raton, FL, 1992. 

\ss
\nind{[F]}  M. Fuchs, {\it The blow-up of $p$-harmonic maps}. Manu. Math. 81 (1993), 
no. 1-2, 89--94.

\ss
\nind{[FMS1]}  A. Freire, S. M\"uller, M. Struwe, {\it  Weak convergence of wave maps from $(1+2)$-dimensional
Minkowski space to Riemannian manifolds.} Invent. Math. 130 (1997), no. 3, 589--617. 

\ss
\nind{[FMS2]}  A. Freire, S. M\"uller, M. Struwe, {\it Weak compactness of wave maps and harmonic maps.}
Ann. Inst. H. Poincaré Anal. Non Linéaire 15 (1998), no. 6, 725--754.

\ss
\nind{[FS]}  C. Fefferman, E. Stein, {\it $H^p$ spaces of several variables}. 
Acta Math. 129 (1972), 137--193.

\ss
\nind{[H1]}  F. H\'elein, {\it Regularite des applications faiblement harmoniques entre une surface
et variete riemannienne}. CRAS, Paris 312 (1991) 591--596. 

\ss
\nind{[H2]}  F. H\'elein, {Harmonic maps, conservation laws and moving frames}. 
Second edition. Cambridge Tracts in Mathematics, 150. Cambridge University Press, Cambridge, 2002.

\ss
\nind{[H]}  N. Hungerbühler, {\it $m$-harmonic flow}. Ann. Scuola Norm. Sup. Pisa Cl. Sci. 
(4) 24 (1997), no. 4, 593--631 (1998).

\ss
\nind{[HL]}  R. Hardt, F. H. Lin, {\it Mappings minimizing the $L\sp p$ norm of the gradient}. 
Comm. Pure Appl. Math. 40 (1987), no. 5, 555--588.

\ss
\nind{[HLM]}  R. Hardt, F. H. Lin, L. Mou, {\it Strong convergence of $p$-harmonic mappings}. 
Progress in partial differential equations: the Metz surveys, 3, 58--64, 
Pitman Res. Notes Math. Ser., 314, Longman Sci. Tech., Harlow, 1994

\ss
\nind{[IM]} T. Iwaniec, G. Martin, {\it Quasiregular mappings in even dimensions}.
Acta Math. 170 (1993), no. 1, 29-81. 

\ss
\nind{[JN]}  F. John, L. Nirenberg, {\it On functions of bounded mean oscillation}. 
Comm. Pure Appl. Math. 14 1961 415--426.

\ss
\nind{[L1]}  P. L. Lions, {\it The concentration-compactness principle in the calculus of
variations: The limit case, I.} Rev. Mat. Iberoamericana 1, no.1 (1985) 145-201.

\ss
\nind{[L2]}  P. L. Lions, {\it The concentration-compactness principle in the calculus of
variations: The limit case, II.} Rev. Mat. Iberoamericana 1, no.2 (1985) 45-121.

\ss
\nind{[L]}  S. Luckhaus, {\it Convergence of minimizers for the $p$-Dirichlet integral}. 
Math. Z. 213 (1993), no. 3, 449--456.

\ss
\nind{[S]}  J. Shatah, {\it Weak solutions and development of singularities of the ${\rm SU}(2)$
$\sigma$-model}. Comm. Pure Appl. Math. 41 (1988), no. 4, 459--469.

\ss
\nind{[SaU]} J. Sacks, K. Uhlenbeck, {\it The existence of minimal immersions of $2$-spheres}.
Ann. of Math. 113 (1981), 1-24.

\ss
\nind{[SU]}  R. Schoen, K. Uhlenbeck, {\it A regularity theory for harmonic maps}. 
J. Differential Geom. 17 (1982), no. 2, 307--335. 

\ss
\nind{[TW]}  T. Toro,  C. Y. Wang, {\it Compactness properties of weakly $p$-harmonic maps
into homogeneous spaces}. Indiana Univ. Math. J. 44 (1995), no. 1, 87--113

\ss
\nind{[U]}  K. Uhlenbeck, {\it Connections with $L^p$-bounds on curvature}. Comm. Math. Phys.
83 (1982), 31-42.
 
\ss
\nind{[W1]}  C. Y. Wang, {\it Bubble phenomena of certain Palais-Smale sequences from surfaces
to general targets}. Houston J. Math. 22 (1996), no. 3, 559--590.

\ss
\nind{[W2]}  C. Y. Wang, {\it Stationary Biharmonic Maps from $\R^n$ into a Riemannian Manifold}.
Comm. Pure Appl. Math., Vol LVII (2004), 0419-0444.

\ss
\nind{[W3]}  C. Y. Wang, {\it Biharmonic maps from $\R^4$ into a Riemannian manifold}. 
Math. Z. (to appear).

\end